\documentclass{amsart}
\usepackage{amssymb, amsmath, amscd}
\usepackage{latexsym}
\usepackage{graphicx}

\newcommand{\al}{\alpha}

\newcommand{\lo}{\longrightarrow}

\newcommand{\BC}{\Bbb C}

\newcommand{\del}{\delta}
\newcommand{\supp}{\operatorname{supp}}

\newcommand{\repn}{representation\,}

\newtheorem{prop}{Proposition}[section]
\newtheorem{theo}{Theorem}[section]
\newtheorem{lemm}{Lemma}[section]

\newtheorem{nota}{Notation}[section]

\begin{document}

\title[restricted algebras]{resticted algebras on inverse semigroups II,
positive definite functions}
\date{}
\author[M. Amini, A.R. Medghalchi]{Massoud Amini, Alireza Medghalchi}
\address{Department of Mathematics, Tarbiat Modarres University,
P.O.Box 14115-175,\linebreak Tehran,Iran,mamini@modares.ac.ir
\linebreak \linebreak \indent Department of Mathematics\\ Teacher
Training University, Tehran, Iran \linebreak
medghalchi@saba.tmu.ac.ir} 
\thanks{This research was supported by Grant 510-2090
of Shahid Behshti University}

\begin{abstract}
The relation between representations and positive definite
functions is a key concept in harmonic analysis on topological
groups. Recently this relation has been studied on topological
groupoids. This is the second in a series of papers in which we
have investigated the concept of "restricted" positive definite
functions and their relation with representations.
\end{abstract}

\maketitle

\section{Introduction.}
In [1] we introduced the concept of restricted representations for
an inverse semigroup $S$ and studied the restricted forms of some
important Banach algebras on $S$. In this paper, we continue our
study by considering the relation between the restricted positive
definite functions and retricted representations. In particular,
we prove restricted versions of the Godement's characterization of
the positive definite functions of finite support (Theorem 2.1).
These results are used in a forthcoming paper to study the
restricted forms of the Fourier and Fourier-Stieltjes algebras on
an inverse semigroup $S$ [2].

All over this paper, $S$ denotes a unital inverse semigroup with
identity $1$. Let us remind that an inverse semigroup $S$ is a
discrete semigroup such that for each $s\in S$ there is a unique
element $s^*\in S$ such that
$$ss^*s=s,\quad s^*ss^*=s^*.$$
The set $E$ of idempotents of $S$  consists of elements the form
$ss^*,\, s\in S$. $E$ is a commutative sub semigroup of $S$. There
is a natural order $\leq$ on $E$ defined by $e\leq f$ if and only
if $ef=e$. A $*$-\repn of $S$ is a pair $\{\pi,\mathcal H_\pi\}$
consisting of a (possibly infinite dimensional) Hilbert space
$\mathcal H_\pi$ and a map $\pi:S\to\mathcal B(\mathcal H_\pi)$
satisfying
$$\pi(xy)=\pi(x)\pi(y),\, \pi(x^*)=\pi(x)^*\quad(x,y\in S),$$
that is a $*$-semigroup homomorphism from $S$ into the inverse
semigroup of partial isometries on $\mathcal H_\pi$. We loosely
refer to $\pi$ as the \repn and it should be understood that there
is always a Hilbert space coming with $\pi$. Let
$\Sigma=\Sigma(S)$ be the family of all $*$-representations $\pi$
of $S$ with $$\|\pi\|:=\sup_{x\in S}\|\pi(x)\|\leq 1.$$ For $1\leq
p<\infty$, $\ell^p(S)$ is the Banach space of all complex valued
functions $f$ on $S$ satisfying
$$\|f\|_p :=\big(\sum_{x\in S} |f(x)|^p\big)^{\frac{1}{p}}<\infty .$$
For $p=\infty$, $\ell^\infty(S)$ consists of those $f$ with
$\|f\|_\infty :=sup_{x\in S} |f(x)|<\infty$. Recall that $\ell
^1(S)$ is a Banach algebra with respect to the product
$$(f*g)(x)=\sum_{st=x} f(s)g(t) \quad
(f,g\in \ell ^1(S)),$$ and $\ell ^2(S)$ is a Hilbert space with
inner product
$$<f,g>=\sum_{x\in S} f(x)\overline{g(x)}\quad (f,g\in\ell^2(S)).$$

Let also put
$$\check{f}(x)=f(x^*),\,  \tilde{f}(x)=\overline{f(x^*)},$$
for each $f\in \ell ^p (S) \quad(1\leq p\leq \infty)$.

As in [1] let us introduce the associated groupoid of an inverse
semigroup $S$. Given $x,y\in S$, the {\bf restricted product} of
$x,y$ is $xy$ if $x^*x=yy^*$, and undefined, otherwise. The set
$S$ with its restricted product forms a groupoid, which is called
the {\bf associated groupoid} of $S$ and we denote it by $S_a$. If
we adjoin a zero element $0$ to this groupoid, and put $0^*=0$, we
get an inverse semigroup $S_r$ with the multiplication rule
$$x\bullet y=\begin{cases}
xy & \text{if}\;\; x^*x=yy^* \\
0 & \text{otherwise}
\end{cases}\quad (x,y\in S\cup\{0\}),$$
which is called the {\bf restricted semigroup} of $S$. A {\bf
restricted representation} $\{\pi,\mathcal H_\pi\}$ of $S$ is a
map $\pi:S\to\mathcal B(\mathcal H_\pi)$ such that
$\pi(x^*)=\pi(x)^*\quad (x\in S)$ and
$$\pi(x)\pi(y)=\begin{cases}
\pi(xy) & \text{if}\;\; x^*x=yy^* \\
0 & \text{otherwise}
\end{cases}\quad (x,y\in S).$$

Let $\Sigma_r=\Sigma_r(S)$ be the family of all restricted
representations $\pi$ of $S$ with $\|\pi\|=\sup_{x\in
S}\|\pi(x)\|\leq 1$. It is not hard to guess that $\Sigma_r(S)$
should be related to $\Sigma(S_r)$. Let $\Sigma_0(S_r)$ be the set
of all $\pi\in \Sigma(S_r)$ with $\pi(0)=0$. Note that
$\Sigma_0(S_r)$ contains all cyclic representations of $S_r$. Now
it is clear that, via a canonical identification,
$\Sigma_r(S)=\Sigma_0(S_r)$. Two basic examples of restricted
representations are the restricted left and right regular
representations $\lambda_r$ and $\rho_r$ of $S$ [1]. For each
$\xi,\eta\in\ell^1(S)$ put
$$(\xi\bullet\eta)(x)=\sum\limits_{x^*x=yy^*}\xi(xy)\eta(y^*)\quad(x\in S),$$
then $(\ell^1(S),\bullet , \tilde\,\,)$ is a semisimple Banach
$*$-algebra [1] which is denoted by $\ell_r^1(S)$ and is called
the restricted semigroup algebra of $S$.

\section{Reduced positive definite functions}

A bounded complex valued function $u:S\lo\BC$ is called {\bf
positive definite} if for all positive integers $n$ and all
$c_1,\dots,c_n\in\BC$, and $x_1,\dots,x_n\in S$, we have
$$\sum_{i=1}^n \sum_{j=1}^n \bar c_ic_j u(x_i^*x_j)\geq 0,$$
and it is called {\bf restricted positive definite} if for all
positive integers $n$ and all $c_1,\dots,c_n\in\BC$, and
$x_1,\dots,x_n\in S$, we have
$$\sum_{i=1}^n \sum_{j=1}^n \bar c_ic_j (\lambda_r(x_i)u)(x_j)\geq 0.$$
We denote the set of all positive definite and restricted positive
definite functions on $S$ by $P(S)$ and $P_r(S)$, respectively.
The two concepts coincide for (discrete) groups.

It is natural to expect a relation between $P_r(S)$ and $P(S_r)$.
Before checking this, note that $S_r$ is hardly ever unital. This
is important, as the positive definite functions in non unital
case should be treated with extra care [5]. Let us take any
inverse semigroup $T$ with possibly no unit. Of course, one can
always adjoin a unit $1$ to $T$ with $1^*=1$ to get a unital
inverse semigroup $T^1=T\cup\{1\}$ (if $T$ happened to have a unit
we put $T^1=T$). However, positive definite functions on $T$ do
not necessarily extend to positive definite functions on $T^1$.
Following [5], we consider the subset $P_e(T)$ of {\bf extendible
positive definite} functions on $T$ which are those $u\in P(T)$
such that $u=\tilde u$ and there exists a constant $c>0$ such that
for all $n\geq 1$ , $x_1,\dots,x_n\in T$ and
$c_1,\dots,c_n\in\mathbb C$,
$$|\sum_{i=1}^{n} c_i u(x_i)|^2\leq \, c\sum_{i=1}^{n}\sum_{j=1}^{n}
\bar c_i c_j u(x_i^*x_j).$$ If $\tau:\ell^\infty(T)\to\ell^1(T)^*$
is the canonical isomorphism, then $\tau$ maps $P_e(T)$ onto the
set of extendible positive bounded linear functionals on
$\ell^1(T)$ (those which are extendible to a positive bounded
linear functional on $\ell^1(T^1)$) and the restriction of $\tau$
to $P_e(T)$ is an isometric affine isomorphism of convex cones [5,
1.1]. Also the linear span $B_e(T)$ of $P_e(T)$ is an algebra [5,
3.4] which coincides with the set of coefficient functions of
$*$-representations of $T$ [5, 3.2]. If $T$ has a zero element,
then so is $T^1$. In this case, we put $P_0(T)=\{u\in P(T):
u(0)=0\}$ and $P_{0,e}(T)=P_0(T)\cap P_e(T)$. To each $u\in
P_e(T)$, there corresponds a cyclic $*$-\repn of $\ell^1(T^1)$
which restricts to a cyclic \repn of $T$ (see the proof of [5,
3.2]). Let $\omega$ be the direct sum of all cyclic
representations of $T$ obtained in this way, then the set of all
coefficient functions of $\omega$ is the linear span of $P_e(T)$
[5, 3.2]. We call $\omega$ the {\bf universal \repn} of $T$.

All these arrangements are for $T=S_r$, as it is an inverse
$0$-semigroup which is not unital unless $S$ is a group. We remind
the reader that our blanket assumption is that $S$ is a unital
inverse semigroup. From now on, we also assume that $S$ has no
zero element (see Example 2.1).

\begin{lemm}
The restriction map $\tau: P_0(S_r)\to P_r(S)$ is an affine
isomorphism of convex cones.
\end{lemm}
{\bf Proof} Let $u\in P(S_r)$. For each $n\geq 1$ ,
$x_1,\dots,x_n\in S_r$ and $c_1,\dots,c_n\in\mathbb C$, we have
$$\sum_{i=1}^{n}\sum_{j=1}^{n} \bar c_i c_j u(x_i^*\bullet x_j)
=\sum_{x_ix_i^*=x_jx_j^*} \bar c_i c_j u(x_i^* x_j)+ u(0)
\big(\sum_{x_ix_i^*\neq x_jx_j^*} \bar c_i c_j\big) ,$$
in particular if $u\in P_0(S_r)$, then
$$\sum_{i=1}^{n}\sum_{j=1}^{n} \bar c_i c_j u(x_i^*\bullet x_j)
=\sum_{i=1}^{n}\sum_{j=1}^{n} \bar c_i c_j (\lambda_r(x_i)u)(x_j),$$
so $\tau$ maps $P_0(S_r)$ into $P_r(S)$.

$\tau$ is clearly an injective affine map. Also if $u\in P_r(S)$
and $v$ is extension by zero of $u$ on $S_r$, then from above
calculation applied to $v$, $v\in P_0(S_r)$ and $\tau(v)=u$, so
$\tau$ is surjective.\qed

It is important to note that the restriction map $\tau$ may fail
to be surjective when $S$ already has a zero element.

\vspace{.3cm} {\bf Example 5.1.} If $S=[0,1]$ with discrete
topology and operations
$$xy=max\{x,y\} , x^*=x\quad (0\leq x,y\leq 1).$$
Then $S$ is a zero inverse semigroup with identity. Here $S_r=S$,
as sets, $P(S)=\{u: u\geq 0, u \, \text{is decreasing}\,\}$ [3],
but the constant function $u=1$ is in $P_r(S)$. This in particular
shows that the map $\tau$ is not necessarily surjective, if $S$
happens to have a zero element. To show that $1\in P_r(S)$ note
that for each $n\geq 1$, each $c_1,\dots, c_n\in \mathbb C$, and
each $x_1,\dots x_n\in S$, if $y_1,\dots, y_k$ are distinct
elements in $\{x_1,\dots, x_n\}$, then for $J_l:=\{j: 1\leq j\leq
n, x_j=y_l\}$, we have $J_i=J_l$, whenever $i\in J_l$, for each
$1\leq i,l\leq k$. Hence

\begin{align*}
\sum_{i,j=1}^{n} \bar c_i c_j \lambda_r(x_i)1(x_j)
&=\sum_{i=1}^{n}\bar c_i \big( \sum_{x_j=x_i}
c_j\big)=\sum_{l=1}^{k}\big(\sum_{i\in J_l} \bar
c_i ( \sum_{j\in J_i} c_j)\big)\\
&=\sum_{l=1}^{k}\big(\sum_{i\in J_l} \bar c_i ( \sum_{j\in J_l}
c_j)\big)=\sum_{l=1}^{k}\big|\sum_{i\in J_l} c_i\big|^2\geq 0.
\end{align*}

\begin{nota} Let $P_{0,e}(S_r)$ be the set of all extendible elements of $P_0(S_r)$.
This is a subcone which is mapped isomorphically onto a subcone
$P_{r,e}(S)$ by $\tau$. The elements of $P_{r,e}(S)$ are called
{\bf extendible restricted positive definite} functions on $S$.
These are exactly those $u\in P_r(S)$ such that $u=\tilde u$ and
there exists a constant $c>0$ such that for all $n\geq 1$ ,
$x_1,\dots,x_n\in S$ and $c_1,\dots,c_n\in\mathbb C$,
$$|\sum_{i=1}^{n} c_i u(x_i)|^2\leq \, c\sum_{x_ix_i^*=x_jx_j^*} \bar c_i c_j u(x_i^*x_j).$$
\end{nota}

\begin{prop} There is an affine isomorphism $\tau$ of convex cones from $P_{r,e}(S)$ onto
$$\ell_r^1(S)_{+}^*\simeq (\mathbb C\delta_0)_{+}^\perp
=\{f\in \ell^\infty(S_r)_{+}: f(0)=0\}=:\ell_0^\infty(S_r)_{+}.$$
\end{prop}
{\bf Proof} The affine isomorphism $\ell_r^1(S)_{+}^*\simeq
\ell_0^\infty(S_r)_{+}$ is just the restriction of the linear
isomorphism of [1, Proposition 4.1] to the corresponding positive
cones. Let us denote this by $\tau_3$. In Notation 2.1 we
presented an affine isomorphism $\tau_2$ from $P_{e,e}(S_r)$ onto
$P_{r,e}(S)$. Finally [5, 1.1], applied to $S_r$, gives an affine
isomorphism from $P_e(S_r)$ onto $\ell^\infty(S_r)_{+}$, whose
restriction is an affine isomorphism $\tau_1$ from $P_{0,e}(S_r)$
onto $\ell_0^\infty(S_r)_{+}$. Now the obvious map $\tau$, which
makes the diagram

\begin{equation*}
\begin{CD}
P_{0,e}(S_r)@>{\tau_1}>>\ell_0^\infty(S_r)_{+}\\
@V{\tau_2}VV    @VV{\tau_3}V\\
P_{r,e}(S)@>>{\tau}>\ell_r^1(S)_{+}^*
\end{CD}
\end{equation*}

commutative, is the desired affine isomorphism.\qed

In [3] the authors developed harmonic analysis on topological
foundation $*$-semigroups (which include all inverse semigroups)
and in particular studied positive definite functions on them. Our
aim in this section is to develop a parallel theory for the
restricted case, and among other results prove the generalization
of the Godement's characterization of positive definite functions
on groups [4] in our restricted context(Theorem 2.1).

For $F,G\subseteq S$, put
$$F\bullet G=\{st: s\in F, t\in G, s^*s=tt^*\}.$$
This is clearly a finite set, when $F$ and $G$ are finite.

\begin{lemm}
If $S$ is an inverse semigroup and $f,g\in\ell^2(S)$, then
$\supp(f\bullet\tilde{g})=(\supp f)\bullet(\supp g)^*$. In
particular, when $f$ and $g$ are of finite supports, then so is
$f\bullet\tilde{g}$.
\end{lemm}
{\bf Proof}
$f\bullet\tilde{g}(x)=\sum\limits_{x^*x=yy^*}f(xy)\overline{g(y)}\neq
 0$ if and only if
$xy\in\supp(f)$, for some $y\in\supp(g)$ with $x^*x=yy^*$. This is
clearly the case if and only if $x=st^*$, for some $s\in supp(f)$
and $t\in supp(g)$ with $s^*s=t^*t$.
 Hence $\supp(f\bullet\tilde{g})=(\supp f)\bullet(\supp g)^*$. \qed

The following lemma follows from the fact that the product
$f\bullet g$ is linear in each variable.

\begin{lemm} (Polarization Identity) For each $f,g\in \ell^2(S)$
\begin{align*}
4f\bullet\tilde{g}&=(f+g)\bullet(f+g)^{\tilde{}}-(f-g)\bullet(f-g)^{\tilde{}}\\
&+i(f+ig)\bullet(f+ig)^{\tilde{}}-i(f-ig)\bullet(f-ig)^{\tilde{}},
\end{align*}
where $i=\sqrt{-1}$. \qed
\end{lemm}

\begin{lemm}
For each $\varphi\in P_{r,f}(S)$, we have
$\tilde\rho_r(\varphi)\geq 0$.
\end{lemm}
{\bf Proof} For each $x,y\in S$
\begin{align*}
<\tilde\rho_r(\varphi)\del_x,\del_y> & =
\sum_z\tilde\rho_r(\varphi)\del_x(z)\overline{\del_y(z)}=
\tilde\rho_r(\varphi)\del_x(y)=\sum_z\varphi(z)\rho_r(z)\del_x(y)\\
&=\sum_{zz^*=y^*y}\varphi(z)\del_x(yz).
\end{align*}
Now if $xx^*=yy^*$ then for $z=y^*x$ we have $zz^*=y^*xx^*y=y^*y$
and conversely $zz^*=y^*y$ and $x=yz$ imply that
$z=zz^*z=y^*yz=y^*x$, and then $x=yy^*x$ and $xz^*=y$, so
$y=xz^*=xx^*y$, that is $yy^*=xx^*yy^*=yy^*xx^*=xx^*$. Hence the
last sum is $\varphi(y^*x)$ if $xx^*=yy^*$, and it is zero,
otherwise. Summing up,
$$
<\tilde\rho_r(\varphi)\delta_{x},\delta_{y}>=(\lambda_r(y)\varphi)(x).
$$
Now for $\xi=\sum_{i=1}^n a_i\delta_{x_i}\in \ell^2_f(S)$, we get
$$<\tilde\rho_r(\varphi)\xi,\xi>=\sum_{i,j=1}^n
a_i
\bar{a}_j<\tilde\rho_r(\varphi)\delta_{x_i},\delta_{x_j}>=\sum_{i,j=1}^n
a_i \bar{a}_j (\lambda_r(x_j)\varphi)(x_i) \geq 0.$$\qed

\begin{lemm} With above notation,
$$\tilde\lambda_r(f)\tilde\rho_r(g)=\tilde\rho_r(g) \tilde\lambda_r(f),$$
for each $f,g\in \ell^1(S)$.
\end{lemm}
{\bf Proof} Given $f,g\in \ell^1(S)$ and $\xi\in \ell^2(S)$, put
$\eta=\tilde\rho_r(g)\xi$ and $\zeta=\tilde\lambda_r(f)\xi$, then
$\eta,\zeta\in\ell^2(S)$ and for each $x\in S$,
\begin{align*}
\tilde\lambda_r(f)\tilde\rho_r(g)\xi(x) & =\sum\limits_{y\in
S}f(y)(\lambda_r(y)\eta)(x)\\
& =\sum\limits_{yy^*=xx^*} f(y)\eta(y^*x)
 =\sum\limits_{yy^*=xx^*}f(y)\sum\limits_{u\in S}
g(u)(\rho_r(u)\xi)(y^*x)\\
& =\sum\limits_{yy^*=xx^*}f(y)\sum\limits_{uu^*=x^*yy^*x}
 g(u)\xi(y^*xu)\\
 & =\sum\limits_{yy^*=xx^*}f(y)\sum\limits_{uu^*=x^*x}
 g(u)\xi(y^*xu),
\end{align*}
and
\begin{align*}
\tilde\rho_r(g)\tilde\lambda_r(f)\xi(x) &= \sum\limits_{u\in S}
g(u)(\rho_r(u)\zeta)(x) =\sum\limits_{uu^*=x^*x} g(u)\zeta(xu)\\
& =\sum\limits_{uu^*=x^*x} g(u)\sum\limits_{y\in S}
f(y)(\lambda_r(y)\xi)(xu)\\
& =\sum\limits_{uu^*=x^*x}g(u)\sum\limits_{yy^*=xuu^*x^*}
 f(y)\xi(y^*xu)\\
& =\sum\limits_{uu^*=x^*x}g(u)\sum\limits_{yy^*=xx^*}
 f(y)\xi(y^*xu),
\end{align*}
which are obviously the same. \qed

\begin{lemm} For each  $\pi\in\Sigma_r(S)$ and each $\xi\in \mathcal H_\pi$,
the coefficient function $u=<\pi(.)\xi,\xi>$ is in $P_{r,e}(S)$.
\end{lemm}
{\bf Proof} For each $n\geq 1$, $c_1,\dots,c_n\in\BC$, and
$x_1,\dots,x_n\in S$, noting that $\pi$ is a restricted \repn, we
have
\begin{align*}
\sum_{i=1}^n \sum_{j=1}^n \bar c_i c_j(\lambda_r(x_i)(u)(x_j) &=
\sum_{x_ix_i^*=x_jx_j^*}
\bar c_i c_j u(x_i^*x_j)\\
&=\sum_{x_ix_i^*=x_jx_j^*}
\bar c_ic_j <\pi(x_i^*x_j)\xi,\xi>\\
&=\sum_{i=1}^n \sum_{j=1}^n\bar c_ic_j<\pi(x_i)^*\pi(x_j)\xi,
\xi>\\
&=\sum_{i=1}^n \sum_{j=1}^n\bar c_ic_j<\pi(x_j)\xi,
\pi(x_i)\xi>\\
&=\|\sum_{i=1}^n c_i \pi(x_i)\xi\|_2^2\geq 0,
\end{align*}
and, regarding $\pi$ as an element of $\Sigma_0(S_r)$ and using
the fact that $u(0)=0$, we have
\begin{align*}
\big|\sum_{k=1}^n c_k u(x_k)\big|^2&=\big|\sum_{k=1}^n  c_k<\pi(x_k)\xi,\xi>\big|^2\\
&\leq \|\xi\|^2\|\sum_{i=1}^n c_i \pi(x_i)\xi\|_2^2\\
&=\|\xi\|^2\sum_{i=1}^n \sum_{j=1}^n
\bar c_i c_j(\lambda_r(x_i)(u)(x_j)\\
&=\|\xi\|^2\sum_{i=1}^n \sum_{j=1}^n \bar c_i c_j(u)(x_i^*\bullet
x_j),
\end{align*}
so $u\in P_{0,e}(S_r)=P_{r,e}(S)$.\qed

The following is proved by R. Godement in the group case [4]. Here
we adapt the proof given in [6].

\begin{theo}
Let $S$ be a unital inverse semigroup. Given
$\varphi\in\ell^\infty(S)$, the following statements are
equivalent.

$(i)$ $\varphi\in P_{r,e}(S)$,

$(ii)$ There is an $\xi\in \ell^2(S)$ such that
$\varphi=\xi\bullet\tilde{\xi}$.

Moreover if $\xi$ is of finite support, then so is $\varphi$.
\end{theo}
{\bf Proof} By above lemma applied to $\pi=\lambda_r$, $(ii)$
implies $(i)$. Also if $\xi\in \ell_f^2(S)$, then by Lemma 2.2,
$\xi\bullet\tilde{\xi}$ is of finite support.

Conversely assume that $\varphi\in P_{r,e}(S)$. Choose an
approximate identity $\{e_\al\}$ for $\ell^1_r(S)$ consisting of
positive, symmetric functions of finite support, as constructed in
[1, Proposition 3.2]. Let $\rho_r$ be the restricted right regular
representation of $S$, then by above lemma
$\tilde\rho_r(\varphi)\geq 0$. Take
$\xi_\al=\tilde\rho_r(\varphi)^{\frac{1}{2}}e_\al\in\ell^2(S)$,
then if $1\in S$ is the identity element, then for each
$\alpha\geq\beta$ we have
\begin{align*}
\|\xi_\al-\xi_\beta\|_2^2 & =
<\tilde\rho_r(\varphi)^{\frac{1}{2}}(e_\al-e_\beta),
\tilde\rho_r(\varphi)^{\frac{1}{2}}(e_\al-e_\beta)>\\
&=<\tilde\rho_r(\varphi)(e_\al-e_\beta), e_\al-e_\beta>
=\varphi\bullet(e_\al-e_\beta) \bullet(e_\al-e_\beta)(1)\\
&\leq\|\varphi\bullet(e_\al-e_\beta)\bullet(e_\al-e_\beta)\|_1
=\|\varphi\bullet(e_\al-e_\beta)\|_1\to
0,
\end{align*}
as $\al,\beta\lo\infty$, where the last equality follows from [1,
Lemma 3.2 $(ii)$]. Hence, there is $\xi\in\ell^2(S)$ such that
$\xi_\al\lo\xi$ in $\ell^2(S)$. Now for each $t\in S$
\begin{align*}
\xi\bullet\tilde{\xi}(t) & = <\lambda_r(t^*)\xi,\xi>
=\lim_\al<\lambda_r(t^*)\tilde\rho_r(\varphi)^{\frac{1}{2}}e_\al,
\tilde\rho_r(\varphi)^{\frac{1}{2}}e_\al>\\
&=\lim_\al<\tilde\rho_r(\varphi)^{\frac{1}{2}}\lambda_r(t^*)
\tilde\rho_r(\varphi)^{\frac{1}{2}}e_\al,e_\al> =\lim_\al<\tilde\rho_r(\varphi)
\lambda_r(t^*)e_\al,e_\al>\\
&=\lim_\al<\lambda_r(t^*)e_\al,\tilde\rho_r(\varphi)e_\al>
=\lim_\al(\bar{\varphi}\bullet(\tilde{e}_\al\bullet
\lambda_r(t^*)e_\al))(1)\\
&=\lim_\al((\bar{\varphi}\bullet e_\al)\bullet
\lambda_r(t^*)e_\al)(1)
=\lim_\al\sum_y(\bar{\varphi}\bullet e_\al)(y)\lambda_r(t^*)e_\al(y^*) \\
&=\lim_\al\sum_y (\bar{\varphi}\bullet e_\al)
(y^*)\lambda_r(t^*)e_\al(y)
=\lim_\al\sum_y (e_\al\bullet\tilde{\varphi})(y) \lambda_r(t^*)e_\al(y) \\
&=\lim_\al<\lambda_r(t^*)e_\al, e_\al\bullet \check{\varphi}>
=\lim_\al e_\al\bullet (e_\al\bullet\check{\varphi})^{\tilde{}}(t) \\
&=\lim_\al e_\al\bullet \bar{\varphi}\bullet
e_\al(t)=\bar{\varphi}(t).
\end{align*}
The last equality follows from the remark after Proposition 3.2 of
[1] and the fact that $|e_\al\bullet \bar{\varphi}\bullet e_\al
(t)-\bar{\varphi}(t)|\leq \|e_\al\bullet\bar{\varphi}\bullet
e_\al-\bar{\varphi}\|_1$.
 Hence
$\varphi=\bar{\xi}\bullet(\bar{\xi})^{\tilde{}}$, as required.
\qed

{\bf acknowledgement.}  The first author would like to thank
hospitality of Professor Mahmood Khoshkam during his stay in
University of Saskatchewan, where the main part of the revision
was done.

\end{document}